\documentclass[10pt,a4paper]{article}
\usepackage{amsmath}
\usepackage{latexsym}
\usepackage{amssymb}
\usepackage{amstext}
\usepackage{fancyhdr}
\usepackage{indentfirst}
\usepackage{graphicx}
\usepackage{epsfig}
\usepackage{pgf}
\usepackage{listings}
\usepackage{cases}
\usepackage{epstopdf}
\usepackage{graphicx}
\date{}

\begin{document}
\title{The normalized Laplacian spectra of the double corona based on $R$-graph}
\author{Ping-Kang Yu, Gui-Xian Tian\footnote{Corresponding author. E-mail address: gxtian@zjnu.cn or guixiantian@163.com.}\\
{\small{\it College of Mathematics, Physics and Information Engineering,}}\\
{\small{\it Zhejiang Normal University, Jinhua, Zhejiang, 321004,
P.R. China}}}\maketitle

\begin{abstract}
For simple graphs $G$, $G_1$ and $G_2$, we denote their double
corona based on $R$-graph by $G^{(R)}\otimes{\{G_1,G_2\}}$. This
paper determines the normalized Laplacian spectrum of
$G^{(R)}\otimes{\{G_1,G_2\}}$ in terms of these of $G$, $G_1$ and
$G_2$ whenever $G$, $G_1$ and $G_2$ are regular. The obtained result
reduces to the normalized Laplacian spectra of the $R$-vertex corona
$G^{(R)}\odot{G_1}$ and $R$-edge corona $G^{(R)}\circleddash{G_2}$
by choosing $G_2$ or $G_1$ as a null-graph, respectively. Finally,
applying the results of the paper, we construct infinitely many
pairs of normalized Laplacian cospectral graphs.
\end{abstract}



\section{Introduction}

Throughout this paper, all graphs considered are finite simple
graphs. Let $G=(V,E)$ be a graph with  vertex set
$V=\{\upsilon_1,\upsilon_2,\dots,\upsilon_n\}$ and edge set $E(G)$.
The adjacency matrix $A(G)$ of $G$ is an $n\times{n}$ matrix whose
$(i,j)$-entry is $1$ if $\upsilon_i$ and $\upsilon_j$ are adjacent
in $G$ and $0$ otherwise. The degree of $\upsilon_i$ in $G$ is
denoted by $d_i=d_{G}(\upsilon_i)$. Let $D(G)$ be the degree
diagonal matrix of G with diagonal entries $d_1, d_2,\dots,d_n$. The
normalized Laplacian matrix $\mathcal{L}(G)$ of $G$ is defined as
$I_n-D(G)^{-1/2}A(G)D(G)^{-1/2}$, where $I_n$ denotes the identity
matrix of order $n$. Denote the characteristic polynomial
$det(xI_n-\mathcal{L}(G))$ of $\mathcal{L}(G)$ by $\phi(G;x)$. Since
$\mathcal{L}(G)$ is a symmetric and positive semi-definite matrix.
Then its eigenvalues, denoted by
$\lambda_1(G),\lambda_2(G),\dots,\lambda_n(G)$, are all real,
non-negative and can be arranged in non-decreasing order
$0=\lambda_1(G)\leq\lambda_2(G)\leq\dots\leq\lambda_n(G)$. The set
of all eigenvalues of $\mathcal{L}(G)$ is called the normalized
Laplacian spectrum of $G$.

The normalized Laplacian matrix $\mathcal{L}(G)$, which is
consistent with the transition probability matrix
$P(G)=D(G)^{-1}A(G)$ in the random walk on $G$ and spectral
geometry\cite{Chung1997}, has attracted people's attention. For
instance, Banerjee and Jost\cite{Banerjee2008} studied how the
normalized Laplacian spectrum is affected by operations such as
motif doubling, graph splitting and joining. Huang and
Li\cite{Huang2015} studied the normalized Laplacian spectrum of some
graph operations, such as subdivision graph, $Q$-graphs, $R$-graphs
and so on. Butler and Grout\cite{Butler2011} constructed many pairs
of non-regular normalized Laplacian cospectral graphs. Chen et
al.\cite{Chen2004} gave an interlacing inequality on the normalized
Laplacian eigenvalues of $G$. Chen and Zhang\cite{Chen2007} obtained
two formulae for the resistance distance and degree-Kirchhoff index
in terms of the normalized Laplacian eigenvalues and eigenvectors of
$G$ and so on. For more review about the normalized Laplacian
spectrum of graphs, readers may refer to \cite{Chung1997}. Recently,
Chen and Liao\cite{Chen2017} determined the normalized Laplacian
spectra of the (edge)corona for two graphs. Furthermore, they also
obtained the degree-Kirchhoff index and the number of spanning trees
of these graphs. In \cite{Das2017}, the normalized Laplacian spectra
of some subdivision-coronas for two regular graphs were computed by
Das and Panigrahi. This paper considers the normalized Laplacian
spectrum of double corona based on $R$-graph. We first recall that
the \emph{R-graph}\cite{Cvetkovic2009} of a graph $G$, denoted by
$G^{(R)}$, is the graph obtained from $G$ by adding a new vertex
corresponding to each edge of $G$ and by joining each new vertex to
the endpoints of the edge corresponding to it. The following graph
operation based on $R$-graph comes from \cite{Barik2017}.
\\
\\
\textbf{Definition 1.1}\cite{Barik2017}. Let $G$ be a connected
graph on $n$ vertices and $m$ edges. Let $G_1$ and $G_2$ be graphs
on $n_1$ and $n_2$ vertices, respectively. The $R$\emph{-graph
double corona} of $G$, $G_1$ and $G_2$, denoted by
$G^{(R)}\otimes{\{G_1,G_2\}}$, is the graph obtained by taking one
copy of $G^{(R)}$, $n$ copies of $G_1$ and $m$ copies of $G_2$, and
then by joining the $i$-th old-vertex of $G^{(R)}$ to every vertex
of the $i$-th copy of $G_1$ and the $j$-th new-vertex of $G^{(R)}$
to every vertex of the $j$-th copy of $G_2$.\\

We remark that here $R$-graph double corona reduces to the
$R$-vertex corona $G^{(R)}\odot{G_1}$ or $R$-edge corona
$G^{(R)}\circleddash{G_2}$ (see \cite{Lan2014} for more information)
whenever we choose $G_2$ or $G_1$ as a null-graph in Definition 1.1,
respectively.

In \cite{Barik2017}, Barik and and Sahoo determined the Laplacian
spectra of $R$-graph double corona for regular graph $G$ and any two
graphs $G_1$ and $G_2$. Song et al.\cite{Song2016} computed the
spectra and Laplacian spectra of double corona based on subdivision
graph. As applications, they determined the number of spanning trees
of the double corona based on subdivision graph and constructed
infinitely many pairs of cospectral(Laplacian cospectral) graphs.
Recently, Lan and Zhou\cite{Lan2014} characterized the spectra,
Laplacian and signless Laplacian spectra of $R$-vertex corona and
$R$-edge corona. At the same time, they also constructed infinitely
many pairs of cospectral, Laplacian cospectral and signless
Laplacian cospectral graphs.

Motivated by the works above , we focus on dertermining the
normalized Laplacian spectrum of $R$-graph double corona
$G^{(R)}\otimes{\{G_1,G_2\}}$ in terms of those of regular graphs $G
,$ $G_1 $ and $G_2$ (see Theorem 2.3). As a special case, we give
the normalized Laplacian spectra of the $R$-vertex corona
$G^{(R)}\odot{G_1}$ and $R$-edge corona $G^{(R)}\circleddash{G_2}$
by choosing $G_2$ or $G_1$ as a null-graph, respectively (see
Corollaries 2.4 and 2.5). Finally, applying these results, we
construct infinitely many pairs of normalized Laplacian cospectral
graphs.

\section{Main results}

In this section , we determine the normalized Laplacian spectrum of
$G^{(R)}\otimes{\{G_1,G_2\}}$ in terms of those of regular graphs
$G$, $G_1$ and $G_2$. To prove our results, we need some
preliminaries. For two matrices $A=(a_{ij})$ and $B=(b_{ij})$ of
same size $m\times{n}$, the Hadamard product $A\circ{B}=(c_{ij})$ of
$A$ and $B$ is a matrix of the same size $m\times{n}$ with entries
$c_{ij}=a_{ij}b_{ij}$ for $i=1,2,\ldots.m$ and $j=1,2,\ldots,n$.
Similarly, the Kronecker product $A\otimes B$ of matrices
$A=(a_{ij})$ of size $m\times{n}$ and $B$ of size $p\times{q}$ is
the $mp\times nq$ partition matrix $a_{ij}B$. It is
proved\cite{Horn} that $AB\otimes{CD}=(A\otimes{C})(B\otimes{D})$,
whenever the products $AB$ and $CD$ exist. Moreover, $(A\otimes
B)^{-1}=A^{-1}\otimes B^{-1}$ for two nonsingular matrices $A$ and
$B$. If $A$ and $B$ are two matrices of order $n$ and $p$
respectively, then $\det(A\otimes B)=(\det A)^{p}(\det B)^{n}$. For
more review about the Kronecker product, see \cite{Horn}.

Throughout this paper, $1_{n}$ denotes the column vector of size $n$
with all the entries equal to one. Let $G$ be a graph on $n$
vertices and $B$ be a matrix of order $n$. For any parameter
$\lambda$, we will use the following notation
\[
\chi_G(B,1_n,\lambda)=1_n^T[\lambda{I_n}-(\mathcal{L}(G)\circ{B})]^{-1}1_n,
\]
which will be simplified as $\chi_G(B)$ for convenience. Remark that
this notation is similar to the $M$-coronal, which is introduced by
Cui and Tian\cite{Cui2012} (also see \cite{McLeman2011}). It is
proved\cite{Cui2012} that if $M$ is a matrix of order $n$ with each
row sum equal to a constant $t$, then
\[
(\lambda I_n-M)^{-1}1_n=\dfrac{1}{\lambda-t}1_n \Rightarrow
1_n^T(\lambda I_n-M)^{-1}1_n=\dfrac{n}{\lambda-t}.
\]

The following Lemmas 2.1 and 2.2 come from
\cite{Cvetkovic2009,Zhang} and
\cite{Das2017}, respectively.\\
\\
\textbf{Lemma 2.1}\cite{Cvetkovic2009,Zhang}. {\it Assume that the
order of all four matrices $M_1$,$M_2$,$M_3$ and $M_4$ satisfy the
rules of operations on matrices. If $M_1$ and $M_4$ are invertible,
then
\[
\begin{array}{l}
 \det \left( {\begin{array}{*{20}c}
   {M_1 } & {M_2 }  \\
   {M_3 } & {M_4 }  \\
\end{array}} \right) = \det M_4  \cdot \det (M_1  - M_2 M_4^{ - 1} M_3 ) \\
\;\;\;\;\;\;\;\;\;\;\;\;\;\;\;\;\;\;\;\;\;\;\;\;\;\;\;\;\;= \det M_1  \cdot \det (M_4  - M_3 M_1^{ - 1} M_2 ). \\
 \end{array}
\]}
\\
\textbf{Lemma 2.2}\cite{Das2017}. {\it If $G$ is an r-regular graph,
then obviously
\[
\mathcal{L}(G)=I_n-\dfrac{1}{r}A(G).
\]}
\\
\textbf{Theorem 2.3.} {\it Let $G$ be an $r$-regular graph with $n$
vertices and $m$ edges. Also let $G_1$ and $G_2$ be $r_1$-regular
and $r_2$-regular with $n_1$ and $n_2$ vertices, respectively.
Assume that $0=\mu_1(G), \mu_2(G), \ldots,\mu_n(G)$;
$0=\eta_1(G_1), \eta_2(G_1), \ldots,\eta_{n_1}(G_1)$ and
$0=\delta_1(G_2), \delta_2(G_2), \ldots,\delta_{n_2}(G_2)$ be the
normalized Laplacian spectra of $G$, $G_1$ and $G_2$, respectively.
Then the normalized Laplacian spectrum of
$G^{(R)}\otimes{\{G_1,G_2\}}$ consists of:
\begin{itemize}
  \item The eigenvalue $\dfrac{1+r_1\eta_j(G_1)}{r_1+1}$ with multiplicity
  $n$, for every eigenvalue $\eta_j(G_1)$ $(j=2,3,\dots,n_1)$ of
  $\mathcal{L}(G_1)$;
  \item The eigenvalue $\dfrac{1+r_2\delta_k(G_2)}{r_2+1}$ with multiplicity
  $m$, for every eigenvalue $\delta_k(G_2)$ $(k=2,3,\dots,n_2)$ of
  $\mathcal{L}(G_2)$;
  \item Four roots of equation
\[\small
\begin{array}{l}
 \left[(x-1)(2+n_2)(xr_2+x-1)-n_2\right]\left[(x-1)(2r+n_1)(xr_1+x-1)+r(1-\mu_i(G))(xr_1+x-1)
  -n_1\right]\\-r(\mu_i(G)-2)(xr_1+x-1)(xr_2+x-1)=0,
\end{array}
\]
  for each eigenvalue $\mu_i(G)$ $(i=1,2,\dots,n)$ of $\mathcal{L}(G)$;
  \item Two roots of equation $(x-1)(2+n_2)(xr_2+x-1)-n_2=0$ with multiplicity $m-n$ if
  $m>n$.\\
\end{itemize}}
\textbf{Proof:} Let $M$ be the vertex-edge incidence matrix of $G$.
Then one has
\begin{equation*}
A(G^{(R)}\otimes{\{G_1 , G_2}\})=
\begin{pmatrix}
A(G) & M & 1_{n_1}^T\otimes{I_n} & 0 \\
M^T & 0 & 0 & 1_{n_2}^T\otimes{I_m} \\
1_{n_1}\otimes{I_n} & 0 & A(G_1)\otimes{I_n} & 0 \\
0 & 1_{n_2}\otimes{I_m} & 0 & A(G_2)\otimes{I_m} \\
\end{pmatrix}
\end{equation*}
and
\begin{equation*}
D(G^{(R)}\otimes{\{G_1 , G_2}\})=
\begin{pmatrix}
(2r+n_1)I_n & 0 & 0 & 0\\
0 & (2+n_2)I_m & 0 & 0\\
0 & 0 & (r_1+1)I_{n_1}\otimes{I_n} &0\\
0 & 0 & 0 &(r_2+1)I_{n_2}\otimes{I_m} \\
\end{pmatrix}.
\end{equation*}
Thus the normalized Laplacian matrix of $G^{(R)}\otimes{\{G_1 ,
G_2}\}$ is
\begin{equation*}
\mathcal{L}(G^{(R)}\otimes{\{G_1 , G_2}\})= \begin{pmatrix}
I_n-\frac{A(G)}{2r+n_1} & -c_1M & -c_21_{n_1}^T\otimes{I_n} & 0\\
-c_1M^T & I_m & 0 & -c_31_{n_2}^T\otimes{I_m}\\
-c_21_{n_1}\otimes{I_n} & 0 & (\mathcal{L}(G_1)\circ{B})\otimes{I_n} & 0\\
0 & -c_31_{n_2}\otimes{I_m} & 0 & (\mathcal{L}(G_2)\circ{C})\otimes{I_m}\\
\end{pmatrix},
\end{equation*}
where
\[
 B=\alpha{J_{n_1}}+(1-\alpha)I_{n_1}\;\text{with}\;\alpha=r_1/(r_1+1);\;
 C=\beta{J_{n_2}}+(1-\beta)I_{n_2} \;\text{with}\; \beta=r_2/(r_2+1)
\]
and
\[
c_1=\dfrac{1}{\sqrt{(2r+n_1)(2+n_2)}},\;c_2=\dfrac{1}{\sqrt{(2r+n_1)(r_1+1)}},\;
c_3=\dfrac{1}{\sqrt{(2+n_2)(r_2+1)}}.
\]
Hence the characteristic polynomial of $G^{(R)}\otimes{\{G_1,G_2\}}$
is $\Phi_{\mathcal{L}(G^{(R)}\otimes{\{G_1,G_2\}})}(x)=\det{B_0}$
where
\[
\begin{array}{l}
B_0=xI-\mathcal{L}(G^{(R)}\otimes{\{G_1,G_2\}})\\
\;\;\;\;\;=\begin{pmatrix}
(x-1)I_n+\dfrac{A(G)}{2r+n_1} & c_1M& c_21_{n_1}^T\otimes{I_n} & 0\\
c_1M^T & (x-1)I_m & 0 & c_31_{n_2}^T\otimes{I_m} & \\
c_21_{n_1}\otimes{I_n} & 0 & (xI_{n_1}-\mathcal{L}(G_1)\circ{B})\otimes{I_n}& 0\\
0 & c_31_{n_2}\otimes{I_m} & 0 & (xI_{n_2}-\mathcal{L}(G_2)\circ{C})\otimes{I_m}\\
\end{pmatrix}.
\end{array}
\]
Denoted by $M_0$ the elementary block matrix below,
\[
M_0=\begin{pmatrix}
I_n & 0 & -c_21_{n_1}^T(xI_{n_1}-\mathcal{L}(G_1)\circ{B})^{-1}\otimes{I_n} &0\\
0& I_m &0 &-c_31_{n_2}^T(xI_{n_2}-\mathcal{L}(G_2)\circ{C})^{-1}\otimes{I_m}\\
0&0& I_{n_1}\otimes{I_n}&0\\
0&0&0& I_{n_2}\otimes{I_m}\\
\end{pmatrix}.
\]
Now, we let $B_1=M_0B_0$. It follows from $\det{M_0}=1$ that
\begin{equation}\label{1}
\begin{array}{l}
\Phi_{\mathcal{L}(G^{(R)}\otimes{\{G_1,G_2\}})}(x)=\det{B_1}\\
\;\;\;\;\;\;\;\;\;\;\;\;\;\;\;\;\;\;\;\;\;\;\;\;\;\;\;\;\;\;\;\;\;=\det(xI_{n_2}-\mathcal{L}(G_2)\circ{C})^{m}
\cdot\det(xI_{n_1}-\mathcal{L}(G_1)\circ{B})^{n}\cdot\det{S},
\end{array}
\end{equation}
where
\[\small
S=\begin{pmatrix}
[x-1-c_2^21_{n_1}^T(xI_{n_1}-\mathcal{L}(G_1)\circ{B})^{-1}1_{n_1}]I_n+\dfrac{A(G)}{2r+n_1} & c_1M\\
c_1M^T&
[x-1-c_3^21_{n_2}^T(xI_{n_2}-\mathcal{L}(G_2)\circ{C})^{-1}1_{n_2}]I_m\\
\end{pmatrix}.
\]
Let
$\chi_{G_1}(B)=1_{n_1}^T(xI_{n_1}-\mathcal{L}(G_1)\circ{B})^{-1}1_{n_1}$
and
$\chi_{G_2}(C)=1_{n_2}^T(xI_{n_2}-\mathcal{L}(G_2)\circ{C})^{-1}1_{n_2}$.
From Lemma 2.1, one obtains
\begin{equation}\label{2}
\det S=\det[(x-1-c_3^2\chi_{G_2}(C))I_m] \cdot\det P,
\end{equation}
where
\begin{equation}\label{3}
\begin{array}{l}
P=(x-1-c_2^2\chi_{G_1}(B))I_n+\dfrac{A(G)}{2r+n_1}-\frac{c_1^2}{x-1-c_3^2\chi_{G_2}(C)}MM^T\\
\;\;\;\;=(x-1-c_2^2\chi_{G_1}(B))I_n+\dfrac{A(G)}{2r+n_1}-\frac{c_1^2}{x-1-c_3^2\chi_{G_2}(C)}(rI_n+A(G)).
\end{array}
\end{equation}

Next we shall compute $\chi_{G_1}(B)$ and $\chi_{G_2}(C)$. From
Lemma 2.2, we get
\begin{equation}\label{4}
\mathcal{L}(G_1)\circ{B}=
I_{n_1}-\dfrac{A(G_1)}{r_1+1}=\dfrac{1}{r_1+1}(I_{n_1}+r_1\mathcal{L}(G_1)).
\end{equation}\\
Similarly,
\begin{equation}\label{5}
\mathcal{L}(G_2)\circ{C}=\dfrac{1}{r_2+1}(I_{n_2}+r_2\mathcal{L}(G_2)).
\end{equation}
Observe that $\mathcal{L}(G_1)1_{n_1}=0$ and
$\mathcal{L}(G_2)1_{n_2}=0$. Then we get
\begin{equation*}
(xI_{n_1}-\mathcal{L}(G_1)\circ{B})1_{n_1}=(x-\dfrac{1}{r_1+1})1_{n_1};
\end{equation*}
$$ (xI_{n_2}-\mathcal{L}(G_2)\circ{C})1_{n_2}=(x-\dfrac{1}{r_2+1})1_{n_2}.$$
Hence,
\begin{equation}\label{6}
\chi_{G_1}(B)=\dfrac{n_1}{x-\dfrac{1}{r_1+1}};\;\;
\chi_{G_2}(C)=\dfrac{n_2}{x-\dfrac{1}{r_2+1}}.
\end{equation}
Now plugging (\ref{6}) into (\ref{3}), again from Lemma 2.2, we have
\begin{equation}\label{7}
\small
 \det P
 =\det \left[ \left(x-1-\dfrac{c_2^2n_1}{x-\dfrac{1}{r_1+1}}+\dfrac{r}{2r+n_1}-\dfrac{c_1^22r}{x-1-c_3^2\dfrac{n_2}{x-\frac{1}{r_2+1}}}\right)I_n
 +\dfrac{c_1^2r\mathcal{L}(G)}{x-1-c_3^2\dfrac{n_2}{x-\frac{1}{r_2+1}}}-\dfrac{r\mathcal{L}(G)}{2r+n_1}\right].
 \end{equation}

Let $\mu_i(G)$, $\eta_j(G_1)$ and $\delta_k(G_2)$ be the eigenvalues
of $\mathcal{L}(G)$, $\mathcal{L}(G_1)$ and $\mathcal{L}(G_2)$,
respectively, for $ i=1,2,\dots,n$; $j=1,2,\dots,n_1$ and
$k=1,2,\dots,n_2$. Then, by (\ref{1}), (\ref{2}), (\ref{7}), along
with (\ref{4}) and (\ref{5}), we obtain
\begin{equation*}\small
\begin{array}{lcl}
\Phi_{\mathcal{L}(G^{(R)}\otimes\{G_1,G_2\})}(x)=\prod\limits_{k=1}^{n_2}\left(x-\dfrac{1+r_2\delta_k(G_2)}
{r_2+1}\right)^m\cdot\prod\limits_{j=1}^{n_1}\left(x-\dfrac{1+r_1\eta_j(G_1)}{r_1+1}\right)^n\cdot\left(x-1-\dfrac{n_2}{(2+n_2)(r_2x+x-1)}\right)^{m-n}\\
\cdot\prod\limits_{i=1}^{n}\left[\left(x-1-\dfrac{n_2}{(2+n_2)(xr_2+x-1)}\right)\left(x-1-\dfrac{n_1}{(2r+n_1)
(xr_1+x-1)}+\dfrac{r(1-\mu_i(G))}{2r+n_1}\right)+\dfrac{r(\mu_i(G)-2)}{(2r+n_1)(2+n_2)}
\right].
\end{array}
\end{equation*}
From the above characteristic polynomial, we have
\begin{itemize}
  \item The eigenvalue $\dfrac{1+r_2\delta_k(G_2)}{r_2+1}$ with multiplicity $m$, for every eigenvalue $\delta_k(G_2)$ $(k=2,3,\dots,n_2)$ of
  $\mathcal{L}(G_2)$;
  \item The eigenvalue $\dfrac{1+r_1\eta_j(G_1)}{r_1+1}$ with multiplicity $n$, for every eigenvalue $\eta_j(G_1)$ $(j=2,3,\dots,n_1)$ of
  $\mathcal{L}(G_1)$;
  \item Four roots of equation
\begin{equation*}
\begin{array}{l}
\left( x-1-\dfrac{n_2}{(2+n_2)(xr_2+x-1)} \right)\left(
x-1-\dfrac{n_1}{(2r+n_1)(xr_1+x-1)}+\dfrac{r(1-\mu_i(G))}{2r+n_1}\right)
\\
+\dfrac{r(\mu_i(G)-2)}{(2r+n_1)(2+n_2)}=0
\end{array}
\end{equation*}
for each eigenvalue $\mu_i(G)$ $(i=1,2,\dots,n)$ of
$\mathcal{L}(G)$;
  \item Two roots of equation
\begin{equation*}
x-1-\dfrac{n_2}{(2+n_2)(xr_2+x-1)}=0
\end{equation*} with multiplicity $(m-n)$ whenever $m>n$.
\end{itemize}
Hence the required result follows. $\Box$\\

Next we consider two special situations of
$G^{(R)}\otimes{\{G_1,G_2\}}$. By choosing $G_2$ as a null-graph, we
can reduce $G^{(R)}\otimes{\{G_1,G_2\}}$ to $R$-vertex corona
$G^{(R)}\odot{G_1}$. Thus, from Theorem 2.3, we obtain\\
\\
\textbf{Corollary 2.4.} {\it Let $G$ be an $r$-regular graph with
$n$ vertices and $m$ edges, $G_1$ be an $r_1$-regular graph with
$n_1$ vertices. Also let $\mu_1(G),\mu_2(G),\dots,\mu_n(G)$ and
$\eta_1(G_1),\eta_2(G_1),\dots,\eta_{n_1}(G_1)$ be the normalized
Laplacian spectra of $G$ and $G_1$, respectively. Then the
normalized Laplacian spectrum of $G^{(R)}\odot{G_1}$ consists of:
\begin{itemize}
  \item The eigenvalue $\dfrac{1+r_1\eta_j(G_1)}{r_1+1}$ with multiplicity $n$ , for every
  eigenvalue $\eta_j(G_1)$ $(j=2,3,\dots,n_1)$ of
  $\mathcal{L}(G_1)$;
  \item The roots of equation
\[\small
  \begin{array}{l}
  2(x-1)\left[(x-1)(2r+n_1)(xr_1+x-1)-n_1+r(1-\mu_i(G))(xr_1+x-1)
  \right]+r(\mu_i(G)-2)(xr_1+x-1)=0
  \end{array}
  \]
for each eigenvalue $\mu_i(G)(i=1,2,\dots,n)$ of $\mathcal{L}(G)$ ,
and
  \item The eigenvalue $1$ with multiplicity $m-n$, if $m>n$.
\end{itemize}}

Instead of choosing $G_2$ as a null-graph, if we choose $G_1$ as a
null-graph, then $G^{(R)}\otimes{\{G_1,G_2\}}$ reduces to $R$-edge
corona $G^{(R)}\circleddash{G_2}$. Thus we arrive at\\
\\
\textbf{Corollary 2.5.} {\it Let $G$ be an $r$-regular graph with
$n$ vertices and $m$ edges, $G_2$ be an $r_2$-regular graph with
$n_2$ vertices. Also let $\mu_1(G),\mu_2(G),\dots,\mu_n(G)$ and
$\delta_1(G_2),\delta_2(G_2),\dots,\delta_{n_2}(G_2)$ be the
normalized Laplacian spectra of $G$ and $G_2$, respectively. Then
the normalized Laplacian spectrum of $G^{(R)}\circleddash{G_2}$
consists of:
\begin{itemize}
  \item The eigenvalue $\dfrac{1+r_2\delta_k(G_2)}{r_2+1}$ with multiplicity $m$, for every eigenvalue $\delta_k(G_2)$ $(k=2,3,\dots,n_2)$ of
  $\mathcal{L}(G_2)$;
  \item Three roots of equation
\[
  \begin{array}{l}
  (2x-1-\mu_i(G))\left[(x-1)(2+n_2)(xr_2+x-1)-n_2
  \right]+(\mu_i(G)-2)(xr_2+x-1)=0
  \end{array}
  \]
for each eigenvalue $\mu_i(G)$ $(i=1,2,\dots,n)$ of
$\mathcal{L}(G)$, and
  \item Two roots of equation
\[
  \begin{array}{l}
 (x-1)(2+n_2)(xr_2+x-1)-n_2=0
   \end{array}
  \]
with multiplicity $m-n$, if $m>n$.
\end{itemize}}
Next we shall present an example to explain our Theorem 2.3.\\
\\
\textbf{Example 2.6.} Let us consider three graphs $G=K_3$ ,
$G_1=P_2$ , and $G_2=P_2$. Then the normalized Laplacian eigenvalues
of $G$ are $(\dfrac{3}{2})^{(2)}$ and $0^{(1)}$, where $a^{ (b)}$
indicates that $a$ is repeated $b$ times. The normalized Laplacian
eigenvalues of $G_1$ and $G_2$ are $2^{(1)}$ and $0^{(1)}$. Applying
Theorem 2.3, the normalized Laplacian spectrum of
$G^{(R)}\otimes{\{G_1,G_2\}}$ consists of:
\begin{itemize}
  \item $(\dfrac{3}{2})^{(3)}$ for the normalized Laplacian
eigenvalue 2 of $G_1$;
  \item $(\dfrac{3}{2})^{(3)}$ for the normalized Laplacian
eigenvalue 2 of $G_2$;
  \item For the normalized Laplacian
eigenvalue $\dfrac{3}{2}$ of $G$, the roots of
$24x^4-76x^3+75x^2-24x+\dfrac{9}{4}=0$ with multiplicity $2$ each,
that is, $(\dfrac{3}{2})^{(2)}$, $(\dfrac{3+\sqrt{3}}{4})^{(2)}$,
$(\dfrac{3-\sqrt{3}}{4})^{(2)}$ and $(\dfrac{1}{6})^{(2)}$.
  \item For the normalized Laplacian
eigenvalue 0 of $G$, the roots of $24x^4-64x^3+48x^2-9x=0$ with
multiplicity $1$ each, that is, $0^{(1)}$,
$(\dfrac{7-\sqrt{13}}{12})^{(1)}$,
$(\dfrac{7+\sqrt{13}}{12})^{(1)}$, $(\dfrac{3}{2})^{(1)}$.
\end{itemize}
On the other hand, according to the computation of $\emph{Matlab}$,
we get directly the normalized Laplacian eigenvalues of
$G^{(R)}\otimes{\{G_1,G_2\}}$ are $0^{(1)}$, $(\dfrac{1}{6})^{(2)}$,
$(\dfrac{3}{2})^{(9)}$, $(\dfrac{3-\sqrt{3}}{4})^{(2)}$,
$(\dfrac{3+\sqrt{3}}{4})^{(2)}$, $(\dfrac{7-\sqrt{13}}{12})^{(1)}$,
$(\dfrac{7+\sqrt{13}}{12})^{(1)}$. This example also shows that
Theorem 2.3 is valid.

Similarly, applying corollary 2.4, the normalized spectrum of
$G^{(R)}\odot{G_1}$ consists of: (1) $(\dfrac{3}{2})^{(3)}$; (2) the
roots of $12x^3-32x^2+24x-\dfrac{9}{2}=0$ with multiplicity $2$
each; (3) the roots of $6x^3-13x^2+6x=0$ with multiplicity $1$ each.
Applying corollary 2.5, the normalized spectrum of
$G^{(R)}\ominus{G_2}$ consists of: (1) $(\dfrac{3}{2})^{(3)}$; (2)
the roots of $16x^3-44x^2+33x-\dfrac{9}{2}=0$ with multiplicity $2$
each; and the roots of $4x^3-8x^2+3x=0$ with multiplicity $1$ each.\\

From above theorem and corollaries, we find that the normalized
spectrum of $R$-graph double corona depends on the degree of
regularities, number of vertices, number of edges and normalized
Laplacian eigenvalues of $G$, $G_1$ and $G_2$. Thus, we can
construct infinitely many pairs of normalized Laplacian cospectral
graphs.\\
\\
\textbf{Lemma 2.7}\cite{Das2017}. {\it Two regular graphs are
normalized Laplacian cospectral if and only if they are cospectral.}\\
\\
\textbf{Theorem 2.8.} {\it If $G$ and $H$ are cospectral regular
graphs (not necessarily distinct), so as to $G_i$ and $H_i$ (for
$i=1,2$) (not necessarily distinct), then
$G^{(R)}\otimes{\{G_1,G_2\}}$ (respectively $G^{(R)}\odot{G_1}$,
$G^{(R)}\ominus{G_2}$) is normalized Laplacian cospectral to
$H^{(R)}\otimes{\{H_1,H_2\}}$ (respectively $H^{(R)}\odot{H_1}$,
$H^{(R)}\ominus{H_2}$).}\\
\\
\textbf{Proof:} From Theorem 2.3 and Lemma 2.7, the result follows.
$\Box$\\
\\
\textbf{Remark 2.9.} The graphs $G$ and $H$, along with $G_i$ and
$H_i$ (for $i=1,2$) are regular in Theorem 2.8, but
$G^{(R)}\otimes{\{G_1,G_2\}}$ and $H^{(R)}\otimes{\{H_1,H_2\}}$ are
non-regular in the general case. Hence we can construct infinitely
many pairs of non-regular normalized Laplacian cospectral graphs by
using double corona operations based on $R$-graphs. In addition, we
remark that the degree Kirchhoff index and the number of spanning
trees of some graph operations have been studied extensively(for
example, see
\cite{Barik2017,Bu2014,Chen2017,Chen2007,Huang2015,Huang2016,Tian2017}).
Our results can also help us to compute the number of spanning trees
and degree Kirchhoff index for $R$-graph double corona operations of
graphs, omitted.\\
\\
\textbf{Acknowledgements} This work was supported in part by
NNSFC(No. 11671053) and the Natural Science Foundation of Zhejiang
Province, China (No. LY15A010011).


\begin{thebibliography}{99}

\bibitem{Banerjee2008} A. Banerjee, J. Jost, On the spectrum of the normalized graph
Laplacian, Linear Algebra Appl. 428 (2008) 3015-3022.

\bibitem{Barik2017} S. Barik, G. Sahoo, On the Laplacian spectra of some variants of corona, Linear Algebra Appl., 512(2017) 32-47.

\bibitem{Bu2014} C.-J. Bu, B. Yan, X.-Q. Zhou, J. Zhou, Resistance distance in
subdivision-vertex join and subdivision-edge join of graphs, Linear
Algebra Appl. 458 (2014) 454-462.

\bibitem{Butler2011} S. Butler, J. Grout, A construction of cospectral graphs for the
normalized Laplacian, Electron. J. Combin. 18 (1) (2011)
$\sharp$P231.


\bibitem{Chen2004} G.T. Chen, G. Davis, F. Hall, Z.S. Li, K. Patel, M. Stewart, An
interlacing result on normalized Laplacians, SIAM J. Discrete Math.
18(2) (2004) 353-361.

\bibitem{Chen2017} H.Y. Chen, L.W. Liao, The normalized Laplacian spectra of
the corona and edge corona of two graphs, Linear Multilinear
Algebra, 65 (2017) 582-592.

\bibitem{Chen2007} H.Y. Chen, F.J. Zhang, Resistance distance and the normalized
Laplacian spectrum, Discrete Appl. Math. 155 (2007) 654-661.

\bibitem{Chung1997} F.R.K. Chung, Spectral Graph Theory, CBMS Regional Conference Series in Mathematics, Amer. Math. Soc., Providence, 1997.

\bibitem{Cui2012} S.-Y. Cui, G.-X. Tian, The spectrum and the signless
Laplacian spectrum of coronae, Linear Algebra Appl., 437 (2012)
1692-1703.

\bibitem{Cvetkovic2009} D. M. Cvetkovi\'{c}, P. Rowinson, H. Simi\'{c}, An
introduction to the Theory of Graph Spectra, Cambridge University
Press, Cambridge, 2009.

\bibitem{Das2017} A. Das, P. Panigrahi, Normalized Laplacian spectrum
of some subdivision-coronas of two regular graphs, Linear
Multilinear Algebra 65 (2017) 962-972.

\bibitem{Horn} R. A. Horn, C. R. Johnson, Topics in matrix analysis, Cambridge University Press, 1991.

\bibitem{Huang2015} J. Huang, S. Li, On the normalized Laplacian spectrum,
degree-Kirchhoff index and spanning trees of graphs, Bull. Aust.
Math. Soc. 91 (2015) 353-367.

\bibitem{Huang2016} J. Huang, S. Li, The normalized Laplacians,
degree-Kirchhoff index and the spanning trees of linear hexagonal
chains, Discrete Appl. Math. 207 (2016) 67-79.

\bibitem{Lan2014} J. Lan, B. Zhou, Spectra of graph operations based on $R$-graph, Linear and Multilinear Algebra, 63(2014) 1401-1422.

\bibitem{McLeman2011} C. McLeman, E. McNicholas, Spectra of coronae, Linear Algebra
Appl. 435 (2011) 998-1007.

\bibitem{Song2016} C.-X. Song, Q.-X. Huang, X.-Y. Huang, Spectra of Subdivision Vertex-edge Corona for Graphs, Advances in Mathematics(China), 45(2016) 37-47.

\bibitem{Tian2017} G.-X. Tian, The asymptotic behavior of (degree-)Kirchhoff indices of iterated total graphs of regular graphs, to
appear in Discrete Appl. Math., (2017).

\bibitem{Zhang} F.-Z. Zhang, The Schur complement and its applications, Springer, 2005.

\end{thebibliography}
\end{document}